\documentclass[11pt]{article}
\usepackage{graphicx}
\usepackage{amsmath}
\usepackage{amssymb}
\usepackage{amsfonts}
\usepackage{enumerate}
\usepackage{latexsym}
\usepackage{epsfig}

\newtheorem{theorem}{Theorem}[section]

\newtheorem{rem}[theorem]{Remark}

\newtheorem{prop}[theorem]{Proposition}

\def\proofbox{\begin{picture}(6.5,6.5)
\put(0,0){\framebox(6.5,6.5){}}\end{picture}}
\newenvironment{proof}{\noindent{\it Proof.\quad}}{\hfill\proofbox}

\begin{document}

\begin{center} {\bf \Large Injective Simplicial Maps of the Arc Complex}

\vspace{0.15in}
\Large{Elmas Irmak and John D. McCarthy}\\
\vspace{0.09in}

\end{center}

\begin{abstract}In this paper, we prove that each injective simplicial map of the arc complex of a compact, connected, orientable surface with nonempty boundary is induced by a homeomorphism of the surface. We deduce, from this result, that the group of automorphisms of  the arc complex is naturally isomorphic to the extended mapping class group of the surface, provided the surface is not a disc, an annulus, a pair of pants, or a torus with one hole. We also show, for each of these special  exceptions, that the group of automorphisms of  the arc complex is naturally isomorphic to the quotient of  the extended mapping class group of the surface by its center.
\end{abstract}




\section{Introduction}

In this paper, $R = R_{g,b}$ will denote a compact, connected, oriented surface of genus $g$ with $b$ boundary components, where $b \geq 1$.
The {\it extended mapping class group}, $\Gamma^*(R)$, of $R$ is the group of isotopy classes of self-homeomorphisms of $R$. The {\it mapping class group}, $\Gamma(R)$, of $R$ is the group of isotopy classes of orientation preserving self-homeomorphisms of $R$. $\Gamma(R)$ is a subgroup of index $2$ in $\Gamma^*(R)$. An arc $A$ on $R$ is called \textit{properly embedded} if
$\partial A \subseteq \partial R$ and $A$ is transversal to
$\partial R$. $A$ is called \textit{nontrivial} (or
\textit{essential}) if $A$ cannot be deformed into $\partial R$ in
such a way that the endpoints of $A$ stay in $\partial R$ during
the deformation. The arc complex $\mathcal{A}(R)$ is the abstract simplicial complex whose simplices are collections of  isotopy classes of properly embedded essential arcs on $R$ which can be represented by disjoint arcs. $\Gamma^*(R)$ acts naturally on $\mathcal{A}(R)$ by simplicial automorphisms of $\mathcal{A}(R)$.

The main results of this paper:

\begin{theorem} Let $R$ be a compact, connected, orientable surface of genus $g$ with
$b \geq 1$ boundary components. If $\lambda : \mathcal{A}(R) \rightarrow \mathcal{A}(R)$ is
an injective simplicial map then $\lambda$ is induced by a homeomorphism $H : R \rightarrow R$
(i.e $\lambda([A]) = [H(A)]$ for every vertex $[A]$ in $\mathcal{A}(R))$.\end{theorem}

\begin{theorem} Let $R$ be a compact, connected, orientable surface of genus $g$ with
$b \geq 1$ boundary components. If $R$ is not a disc, an annulus, a pair of pants or a torus with one hole, then $Aut(\mathcal{A}(R)) \cong \Gamma^*(R)$. For each of these special cases $Aut(\mathcal{A}(R)) \cong \Gamma^*(R)/Z(\Gamma^*(R))$.  \end{theorem}

\protect\nopagebreak\noindent\rule{1.5in}{.01in}{\vspace{0.005in}}

\small{The authors thank the Max Planck Institute for Mathematics, Bonn for the excellent conditions provided for their stay at this institution, during which these results were obtained.}

Extended mapping class group was viewed as the automorphism group
of the curve complex on orientable surfaces by Ivanov to get information on the algebraic structure
of the mapping class groups. Ivanov proved that the automorphism group of
the curve complex is isomorphic to the extended mapping class group for connected orientable surfaces of genus at least 2 in \cite{ivanov1}. In his paper, he also proved that the automorphisms of the arc complex which are induced by automorphisms of the curve complex are induced by homeomorphisms of the surface. As an application he proved that isomorphisms between any two finite index subgroups are geometric if the genus is at least 2. Ivanov's results were proven by Korkmaz in \cite{K} and Luo in \cite{L} for lower genus cases.

After Ivanov's work, extended mapping class group was viewed as the automorphism group
of various geometric objects on surfaces. Some of these objects include Schaller's complex (see \cite{S} by Schaller), the
complex of pants decompositions (see \cite{M} by Margalit), the complex of nonseparating curves (see \cite{Ir3} by Irmak), the complex of separating curves (see \cite{BM1} by Brendle-Margalit, and \cite{MV} by McCarthy-Vautaw), the complex of Torelli geometry (see \cite{FIv} by Farb-Ivanov), the Hatcher-Thurston complex (see \cite{IrK} by Irmak-Korkmaz). As applications, Farb-Ivanov proved that the automorphism group of the Torelli subgroup is isomorphic to the mapping class group in \cite{FIv}, and McCarthy-Vautaw extended this result to $g \geq 3$ in \cite{MV}.

Some similar results on simplicial maps and the applications for orientable surfaces are as follows: Irmak proved that superinjective simplicial maps of the curve complex are induced by homeomorphisms of the surface to classify injective homomorphisms from finite
index subgroups of the mapping class group to the whole group (they are geometric except for closed genus two surface) for genus at least two in \cite{Ir1}, \cite{Ir2}, \cite{Ir3}. Behrstock-Margalit and Bell-Margalit proved these results for lower genus cases in \cite{BhM} and in \cite{BeM}. Brendle-Margalit proved that superinjective simplicial maps of separating curve complex are induced by homeomorphisms, to prove that an injection from a finite index subgroup of
$K$ to the Torelli group, where $K$ is the subgroup of mapping class group
generated by Dehn twists about separating curves, is induced by a
homeomorphism in \cite{BM1}, \cite{BM2}. Shackleton proved that
injective simplicial maps of the curve complex are induced by homeomorphisms in \cite{Sh} (he also considers maps between different surfaces), and he obtained strong local co-Hopfian results for mapping class groups. Bell-Margalit proved that superinjective simplicial maps of the curve complex are onto. For nonorientable odd genus surfaces Atalan-Ozan proved that the automorphism group of the curve complex is isomorphic to the mapping class group if $g + r \geq 6$ in \cite{A}.

After this paper was written, Irmak, the first author of this paper, proved similar results about the arc complex on nonorientable surfaces in \cite{Ir4}. The main results in \cite{Ir4}: Let $N$ be a compact, connected, nonorientable surface of genus $g$ with
$r \geq 1$ boundary components. If $\lambda : \mathcal{A}(N) \rightarrow \mathcal{A}(N)$ is
an injective simplicial map then $\lambda$ is induced by a homeomorphism $h : N \rightarrow N$, and
$Aut(\mathcal{A}(N)) \cong \Gamma(N) /Z(\Gamma(N))$.

\section{Mapping Class Groups and Complex of Arcs}

In this section we will prove our main results for $(g, b) = (0,1)$, $(g, b)=(0,2)$, $(g, b)=(0,3)$, and $(g, b)=(1,1)$. In the next section we will give a general argument for the proof of the remaining cases. Unless otherwise indicated, all arcs will be assumed to be essential arcs on $R$.
We will denote arcs by capital letters and their isotopy classes by the corresponding lower case letters (e.g. $A$ and $a =[A] \in \mathcal{A}(R)$).

\begin{theorem}
\label{mainthm1} Let $R$ be a compact, connected, orientable surface of genus $g$ with
$b$ boundary components. If $(g, b) \in \{(0,1), (0,2), (0,3),$ $ (1,1)\}$
then $Aut(\mathcal{A}(R)) \cong \Gamma^*(R) /Z(\Gamma^*(R))$. \end{theorem}

\begin{proof} Case (i): Suppose that $(g, b) = (0,1)$. $R$ is a disc and no arc on $R$ is essential. Hence, $\mathcal{A}(R) = \emptyset$; every injective simplicial map $\mathcal{A}(R) \rightarrow \mathcal{A}(R)$ is an automorphism of $\mathcal{A}(R)$; and $Aut(\mathcal{A}(R))$ is a trivial group. $\Gamma^*(R)$ is a cyclic group of order two. It follows that  $Z(\Gamma^*(R)) = \Gamma^*(R)$ and, hence, $\Gamma^*(R)/Z(\Gamma^*(R))$ is also a trivial group. Hence, $\Gamma^*(R)/Z(\Gamma^*(R))$ is isomorphic to $Aut(\mathcal{A}(R))$.

Case (ii): Suppose that $(g, b) = (0,2)$. $R$ is an annulus, $\mathcal{A}(R)$ consists of a single vertex and $Aut(\mathcal{A}(R))$ is a trivial group. The action of $\Gamma^*(R)$ on $\pi_0(\partial R)$ yields a short exact sequence:
\begin{equation} 1 \rightarrow Z_2 \rightarrow \Gamma^*(R) \rightarrow \Sigma(\pi_0(\partial R)) \rightarrow 1
\label{equation:annulus} \end{equation}

\noindent where $\Sigma(\pi_0(\partial R)) \cong \Sigma_2$ is the group of permutations of $\pi_0(\partial R)$, $\Gamma^*(R) \rightarrow \Sigma(\pi_0(\partial R))$ is the corresponding representation, and the kernel $Z_2$ of $\Gamma^*(R) \rightarrow \Sigma(\pi_0(\partial R))$ is the cyclic group of order $2$ generated by the isotopy class of any orientation reversing involution of $R$ which preserves each component of $\partial R$. The natural representation $\Gamma^*(R) \rightarrow \Sigma(\pi_0(\partial R))$ restricts to an isomorphism $\Gamma(R) \rightarrow \Sigma(\pi_0(\partial R))$. It follows that the above exact sequence (\ref{equation:annulus}) is a split short exact sequence, and $\Gamma^*(R)$ is isomorphic to $\mathbb{Z}_2 \oplus \mathbb{Z}_2$. This implies that the center $Z(\Gamma^*(R))$ of $\Gamma^*(R)$ is equal to $\Gamma^*(R)$; $\Gamma^*(R)/Z(\Gamma^*(R))$ is a trivial group; and, hence, $\Gamma^*(R)/Z(\Gamma^*(R))$ is isomorphic to $Aut(\mathcal{A}(R))$.

\begin{figure}[!hbp]
  \begin{center}
  \scalebox{0.40}{\includegraphics{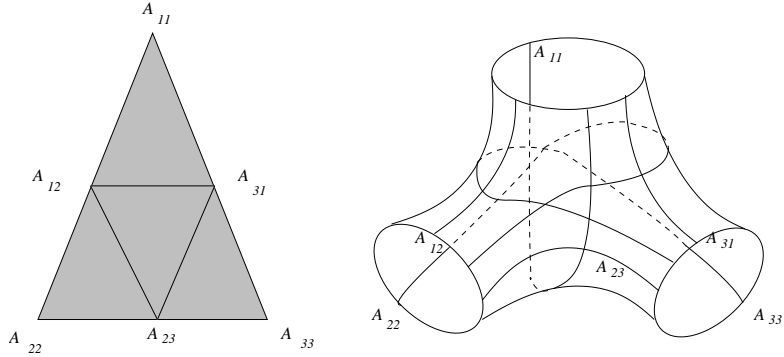}}
  \end{center}
   \caption{The arc complex of a sphere with three holes and
   arcs representing its six vertices}
  \label{fig:A03}\end{figure}

Case (iii): Suppose that $(g, b) = (0,3)$. $R$ is a pair of pants (i.e. a sphere with three holes), then there are exactly six isotopy classes of essential arcs on $R$ and $\mathcal{A}(R)$ is a two-complex represented by a regular tessellation of a triangle by four triangles as illustrated in Figure \ref{fig:A03}. $Aut(\mathcal{A}(R))$ is isomorphic to the symmetric group $\Sigma_3$ on three letters. Indeed, $Aut(\mathcal{A}(R))$ is naturally isomorphic to the group of permutations $\Sigma(\pi_0(\partial R))$ of the set of components $\pi_0(\partial R)$ of $\partial R$. The action of $\Gamma^*(R)$ on $\pi_0(\partial R)$ yields a short exact sequence:
\begin{equation} 1 \rightarrow Z_2 \rightarrow \Gamma^*(R) \rightarrow \Sigma(\pi_0(\partial R)) \rightarrow 1
\label{equation:pants} \end{equation}

\begin{figure}[!hbp]
  \begin{center}
  \scalebox{0.45}{\includegraphics{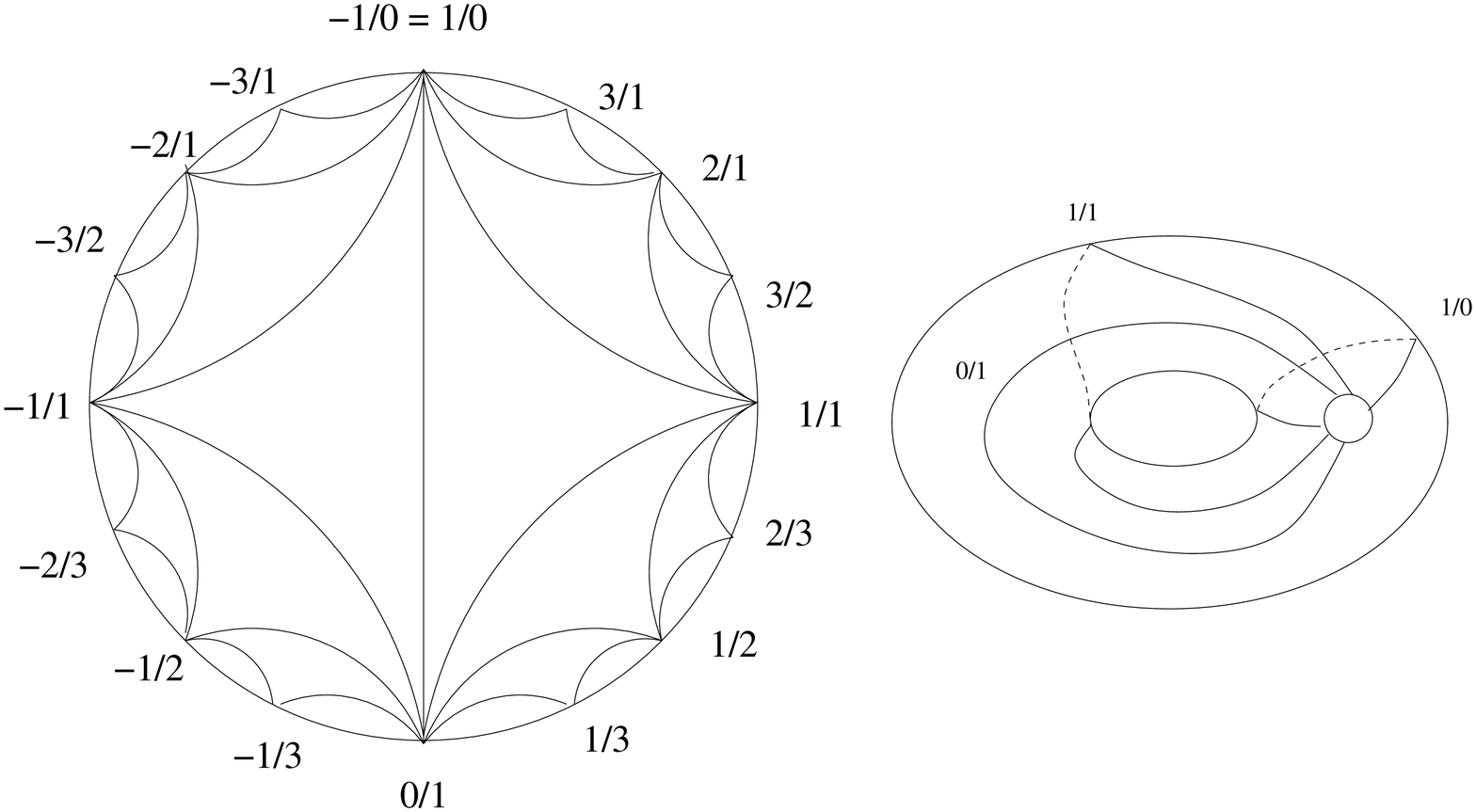}}
  \end{center}
   \caption{The arc complex of a torus with one hole and
   three arcs representing one of its triangles}
  \label{fig:A11}\end{figure}

\noindent where $\Sigma(\pi_0(\partial R)) \cong \Sigma_3$ is the group of permutations of $\pi_0(\partial R)$, $\Gamma^*(R) \rightarrow \Sigma(\pi_0(\partial R))$ is the corresponding representation, and the kernel $Z_2$ of $\Gamma^*(R) \rightarrow \Sigma(\pi_0(\partial R))$ is the cyclic group of order $2$ generated by the isotopy class of any orientation reversing involution of $R$ which preserves each component of $\partial R$. The natural representation $\Gamma^*(R) \rightarrow \Sigma(\pi_0(\partial R))$ restricts to an isomorphism $\Gamma(R) \rightarrow \Sigma(\pi_0(\partial R))$. It follows that the above exact sequence (\ref{equation:pants}) is a split short exact sequence. Since $\Sigma_3$ has trivial center, it follows that $Z_2$ is equal to the center $Z(\Gamma^*(R))$ of $\Gamma^*(R)$; $\Gamma^*(R)/Z(\Gamma^*(R))$ is also naturally isomorphic to $\Sigma(\pi_0(\partial R))$. Hence,
$\Gamma^*(R)/Z(\Gamma^*(R))$ is naturally isomorphic to $Aut(\mathcal{A}(R))$.

Case (iv): Suppose that $(g, b) = (1,1)$. $R$ is a torus with one hole. $\mathcal{A}(R)$ is represented by the decomposition of the hyperbolic plane $\mathbb{H}$ into ideal triangles by the familiar {\it Farey graph}, $\mathcal{F}$ (see Figure 2).

More precisely, let $S$ be the torus obtained by attaching a disc $D$ to $\partial R$ and $P$ be a point in the interior of $D$. Choose an identification of $(S,P)$ with the standard torus, $(S^1 \times S^1,(1,1))$. Then the isotopy classes of arcs on $R$ correspond naturally to the rational points on the circle at infinity $S_\infty = \mathbb{R}^* = \mathbb{R} \cup \infty$, where the arc $A$ on $R$ corresponds to the rational point $p/q$ if and only if the extension of the arc $A$ on $R$ to a closed curve on $S$ by ``coning off'' the endpoints of $A$ in $\partial D$ to the ``center'' $P$ of $D$ represents $\pm (p,q) \in \mathbb{Z} \oplus \mathbb{Z} = \pi_1(S^1 \times S^1,(1,1))$.

The ideal triangles of the decomposition of the hyperbolic plane $\mathbb{H}$ by the {\it Farey graph} correspond to ideal triangulations of $(S,P)$, which correspond to maximal simplices of $\mathcal{A}(R)$. As is well-known, $(\Gamma^*(R),Z(\Gamma^*(R)) \cong (GL(2,\mathbb{Z}),\{\pm I\})$, and, hence, $\Gamma^*(R)/Z(\Gamma^*(R)) \cong PGL(2,\mathbb{Z}) \cong Aut(\mathcal{F}) \cong $ $ Aut(\mathcal{A}(R))$.
\end{proof}

\begin{theorem}
\label{?} Let $R$ be a compact, connected, orientable surface of genus $g$ with
$b$ boundary components. Suppose that $(g, b) \in \{(0,1), (0,2),$ $ (0,3), (1,1)\}$. If $\lambda : \mathcal{A}(R) \rightarrow \mathcal{A}(R)$ is
an injective simplicial map then $\lambda$ is induced by a homeomorphism $H : R \rightarrow R$. \end{theorem}

\begin{proof} For the cases $(g, b)=(0,1), (g, b)=(0,2)$ and $(g, b)=(0,3)$ the proof follows from Theorem 2.1 as for all of these cases the arc complex has finitely many vertices so every injective simplicial map is an automorphism.

When $(g, b)=(1,1)$, we will prove that every injection is onto. Let $w$ be a vertex of $\mathcal{A}(R)$.
Let $\Delta$ be a top dimensional simplex. Since $\lambda$ is injective, we see that $\lambda(\Delta)$ corresponds to a top dimensional simplex, $\Delta'$ in $\mathcal{A}(R)$. If $w$ is a vertex of $\Delta'$, then $w$ is in the image. Suppose that $w$ is not a vertex of $\Delta'$.
Take a top dimensional simplex $\Delta''$ containing $w$. It is easy to see that there exists a chain
$\Delta' = \Delta_0', \Delta_1', \cdots, \Delta_n' = \Delta''$ of top dimensional simplices in $\mathcal{A}(R)$, connecting $\Delta'$ to $\Delta''$ in such a way that the consecutive simplices
$\Delta_i', \Delta_{i+1}'$ have exactly one common face of codimension 1. Let $x'$ be a vertex of $\Delta'$  which is not in $\Delta_1'$. Let $x, y, z$ be the vertices of $\Delta$ such that $\lambda(x) = x'$. There exists a unique top dimensional simplex containing $y, z$ and not $x$, call it $\Delta_1$. It is easy to see that $\lambda(\Delta_1)= \Delta_1'$, and so every vertex in $\Delta_1'$ is in the image of $\lambda$. By an inductive argument, using the above sequence we see that $w$ is in the image of $\lambda$. The result of the theorem now follows from Theorem 2.1.
\end{proof}

\section{Triangulations}

In this section, we assume that $(g, b) \neq (0,1), (g, b) \neq (0,2)$, $(g, b) \neq (0,3)$ and $(g, b) \neq (1,1)$.
Let $\lambda : \mathcal{A}(R) \rightarrow \mathcal{A}(R)$ be an injective simplicial map. We will prove some properties of $\lambda$. First we give some definitions.

Let $T$ be a set of pairwise disjoint nonisotopic arcs on $R$. $T$ is called a {\it triangulation of} $R$ if each component $\Delta$ of the surface $R_T$, obtained from $R$ by cutting $R$ along $T$, is a disc with boundary $\partial \Delta$ equal to a union of arcs, $A, B, C, D, E$, and $F$, where $A$, $B$, and $C$ correspond to elements of $T$ and $D$, $E$, and $F$ correspond
to arcs or circles in $\partial R$. $\Delta$ is called a {\it triangle of} $T$, and $A, B, C$ are called sides of $\Delta$. If $A$, $B$, and $C$
correspond to distinct elements of $T$, then $\Delta$ is called an {\it an embedded triangle of} $T$. Otherwise, $\Delta$ is called {\it a non-embedded triangle of} $T$. The phrase {\it triangle of} $T$ will also be used to refer to the image of any component $\Delta$ of $R_T$, under the natural quotient map $q : R_T \rightarrow R$. The images of $A, B$ and $C$ will also be called as sides of the image triangle. Two distinct triangles of a triangulation $T$ are called {\it adjacent} w.r.t. $T$ if they have a common side.

Let $T$ be a triangulation of $R$. Let $[T]$ be the set of isotopy classes of elements of $T$. Note that $[T]$ is a maximal simplex of $\mathcal{A}(R)$. Every maximal simplex $\sigma$ of $\mathcal{A}(R)$ is equal to $[T]$ for some triangulation $T$ of $R$. So, $\lambda([T]) = [T']$ for some triangulation $T'$ of $R$ and $\lambda$ restricts to a correspondence $\lambda| : [T] \rightarrow [T']$
on the isotopy classes. Note that the triangulation $T'$ of $R$ is well defined up to isotopy on $R$. By using Euler
characteristic arguments we see that the number of arcs in a triangulation is $6g+3b-6$, and the number of triangles in a triangulation is $4g+2b-4$.

Let $a$ and $b$ be isotopy classes of properly embedded essential arcs on $R$. The geometric intersection number $i(a,b)$ of $a$ and $b$ is the minimum number of points in $A \cap B$ where $A$ and $B$ are arcs on $R$ which represent $a$ and $b$.

\begin{prop} If $a$ and $b$ are two vertices of the complex of arcs $\mathcal{A}(R)$ such that $i(a,b)=1$, then $i(\lambda(a), \lambda(b)) = 1$.\label{prop:invgeomint} \end{prop}

\begin{proof} Let $A$ and $B$ be representatives of $a$ and $b$ intersecting once. Note that we may complete $A$ to a triangulation $T_1$ of $R$ such that $(T_1 \setminus \{A\}) \cup \{B\}$ is also a triangulation of $R$. Let $T_2 = (T_1 \setminus \{A\}) \cup \{B\}$.
Let $\sigma_i$ be the simplex of $\mathcal{A}(R)$ corresponding to the triangulation $T_i$ of $R$, $i = 1,2$, and $\sigma_0 = \sigma_1 \cap \sigma_2$. Note that $\sigma_2 \setminus \{b\} = \sigma_0 = \sigma_1 \setminus \{a\}$, and $\sigma_2$ is obtained from $\sigma_1$ by replacing $a$ with $b$ (an elementary move).

Let $\sigma'_i = \lambda(\sigma_i)$, $i = 0,1,2$, $a' = \lambda(a)$, and $b' = \lambda(b)$. Since $\lambda$ is injective there exists a triangulation $T'_i$ corresponding to $\sigma'_i$, $i = 1, 2$.
Since $i(a,b) \neq 0$, there does not exist a simplex of $\mathcal{A}(R)$ having both $a$ and $b$ as vertices. Since $a \in \sigma_1$ and $b \in \sigma_2$, it follows that $\sigma_1 \neq \sigma_2$. Since $\lambda : \mathcal{A}(R) \rightarrow \mathcal{A}(R)$ is an injective simplicial map, it follows that $\sigma'_1 \neq \sigma'_2$.
Let $A'$ be the representative of $a'$ in $T'_1$.
Since $\sigma_2 \setminus \{b\} = \sigma_0 = \sigma_1 \setminus \{a\}$ and $\lambda : \mathcal{A}(R) \rightarrow \mathcal{A}(R)$ is an injective simplicial map, $\sigma'_2 \setminus \{b'\} = \sigma'_0 = \sigma'_1 \setminus \{a'\}$.

Note that we may choose a representative $B'$ of $b'$ such that $B'$ is disjoint from and not isotopic to each element of $T'_1 \setminus \{A'\}$. Let $T'_2 = (T'_1 \setminus \{A'\}) \cup \{B'\}$. Then $T'_2$ is a triangulation of $R$ and $\sigma'_2$ is the simplex of $\mathcal{A}(R)$ corresponding to $T'_2$. Since $\sigma'_1$ and $\sigma'_2$ are distinct maximal simplices of $\mathcal{A}(R)$ containing $\sigma'_0$, $\sigma'_2$ is obtained from $\sigma'_1$ by an elementary move replacing $a'$ with $b'$, we see that $i(a',b') = 1$, completing the proof.
\end{proof}

\begin{prop} Let $\Delta$ be an embedded triangle on $R$ with sides corresponding to $A$, $B$ and $C$. Then there exists a triangulation $T$ on $R$ containing $\{A,B,C\}$ such that the unique triangles $\Delta_A$, $\Delta_B$ and $\Delta_C$ of $T$ on $R$ which are different from $\Delta$ and have, respectively, a side corresponding to $A$, a side corresponding to $B$, and a side corresponding to $C$, are distinct triangles of $T$ on $R$.
\label{prop:fourtriangles} \end{prop}

\begin{proof} Since $\Delta$ is an embedded triangle on $R$ with sides corresponding to $A$, $B$, and $C$, these are nonisotopic essential properly embedded arcs on $R$. There exists a triangulation $T$ of $R$ such that $\{A,B,C\}$ is contained in $T$. Since $\Delta$ is an embedded triangle on $R$, there exist unique triangles $\Delta_A$, $\Delta_B$, and $\Delta_C$ of $T$ on $R$ which are different from $\Delta$ and have, respectively, a side corresponding to $A$, a side corresponding to $B$, and a side corresponding to $C$. Since $R$ is not a pair of pants, $\Delta_A$, $\Delta_B$, and $\Delta_C$ are not the same triangle on $R$.

Suppose that $\Delta_A = \Delta_B$. Then $\Delta_A$ has a side corresponding to $A$ and another side corresponding to $B$.

Suppose that $\Delta_A$ is a non-embedded triangle on $R$. Then, since $\Delta_A$ has sides corresponding to $A$ and $B$, either $\Delta_A$ is the unique triangle of $T$ on $R$ having a side corresponding to $A$ or $\Delta_A$ is the unique triangle of $T$ on $R$ having a side corresponding to $B$, which is a contradiction, as $\Delta$ is a triangle different from $\Delta_A$ having a side corresponding to $A$ and a side corresponding to $B$. It follows that $\Delta_A$ is an embedded triangle on $R$.

$\Delta_A$ has a side corresponding to an element $D$ of $T$, where $D$ is not equal to $A$ or $B$. Suppose that $D = C$. Then $\Delta_A$ is a triangle of $T$ on $R$ different from $\Delta$ having a side corresponding to $C$. In other words, $\Delta_A = \Delta_C$ and, hence, $\Delta_A$, $\Delta_B$, and $\Delta_C$ are the same triangle on $R$, which is a contradiction. Hence, $D$ is not equal to $C$.

\begin{figure}[!hbp]
  \begin{center}
  \scalebox{0.35}{\includegraphics{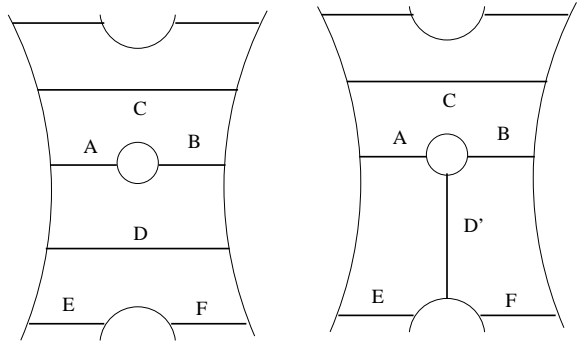}}
  \end{center}
   \caption{Obtaining four triangles by one elementary move}
  \label{fig:fourtriangles1}\end{figure}

Since $\Delta_A$ is an embedded triangle of $T$ on $R$ with sides corresponding to $A$, $B$, and $D$, there exists a unique triangle $\Delta_D$ of $T$ on $R$ which is different from $\Delta_A$ and has a side corresponding to $D$. Note that there is exactly one side of $\Delta_D$ corresponding to $D$. Suppose that the other two sides of $\Delta_D$ correspond to elements $E$ and $F$ of $T$.

Suppose, on the one hand, that $C$ is not equal to $E$ or $F$. Since the sides of $\Delta_D$ correspond to $D$, $E$, and $F$, none of which are equal to $C$, $\Delta_D$ has no side corresponding to $C$. Since $\Delta_C$ has a side corresponding to $C$, it follows that $\Delta_C$ and $\Delta_D$ are distinct triangles of $T$ on $R$. Since $\Delta_A$ and $\Delta_D$ are distinct triangles of $T$ on $R$ having a side corresponding to $D$, there is a quadrilateral $\Omega$ on $R$ with sides corresponding to $A$, $B$, $E$, and $F$, and diagonal $D$ as shown in the first part of Figure 3. Let $D'$ be a diagonal of $Q$ such that $\{D,D'\}$ is a pair of diagonals of $\Omega$ as shown in the second part of Figure 3. Let $T' = (T \setminus \{D\}) \cup \{D'\}$ be the triangulation on $R$ which is obtained from the triangulation $T$ on $R$ by an elementary move replacing $D$ with $D'$. It follows that the unique triangles $\Delta'_A$, $\Delta'_B$, and $\Delta'_C$ of $T'$ on $R$ which are distinct from the triangle $\Delta$ of $T'$ on $R$ and have, respectively, a side corresponding to $A$, a side corresponding to $B$, and a side corresponding to $C$ are distinct triangles of $T'$ on $R$ (see Figure \ref{fig:fourtriangles1}).

Suppose, on the other hand, that $C$ is equal to either $E$ or $F$. We may assume that $C = E$. It follows, by arguments similar to those given above, that $\Delta_D$ is an embedded triangle of $T$ on $R$ with sides corresponding to $C$, $D$, and $F$, where $F$ is some element of $T$ which is not equal to $A$, $B$, $C$, or $D$.

Since $F$ is a side of the embedded triangle $\Delta_D$ of $T$ on $R$, there is a unique triangle $\Delta_F$ of $T$ on $R$ which is distinct from $\Delta_D$ and has a side corresponding to $F$. By arguments similar to those given above, there is exactly one side of $\Delta_F$ which corresponds to $F$. Let the other two sides of $\Delta_F$ correspond to elements $G$ and $H$ of $T$.

\begin{figure}[!hbp]
  \begin{center}
  \scalebox{0.35}{\includegraphics{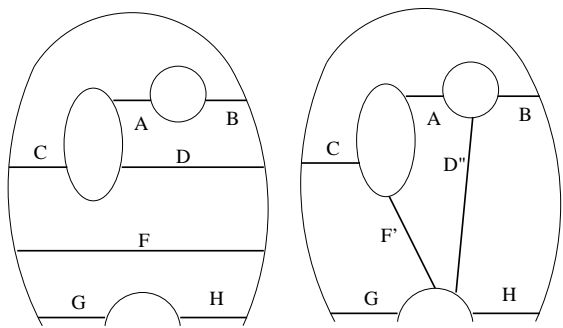}}
  \end{center}
   \caption{Obtaining four triangles by two elementary moves}
  \label{fig:fourtriangles2}\end{figure}

Let $T'$ be the triangulation obtained from $T$ by an elementary move replacing the element $F$ of $T$ by an element $F'$ of $T$. Then let $T''$ be the triangulation obtained from $T'$ by an elementary move replacing the element $D$ of $T'$ by an element $D''$ of $T''$. It follows that the unique triangles $\Delta''_A$, $\Delta''_B$, and $\Delta''_C$ of $T''$ on $R$ which are distinct from the triangle $\Delta$ of $T''$ on $R$ and have, respectively, a side corresponding to $A$, a side corresponding to $B$, and a side corresponding to $C$, are distinct triangles of $T''$ on $R$ (see Figure \ref{fig:fourtriangles2}).

This shows, in any case, that there exists a triangulation of $R$ with the desired properties, completing the proof.\end{proof}\\

Let $\{a,b,c\}$ be a $2$-simplex of $\mathcal{A}(R)$. We say that $\{a,b,c\}$ corresponds to an embedded triangle on $R$ if there exists an embedded triangle $\Delta$ on $R$ with sides corresponding to $A$, $B$, and $C$ representing $a$, $b$, and $c$.\label{defn:simpembedded}

\begin{figure}[htb]
\begin{center}
\epsfxsize=10cm
\epsfbox{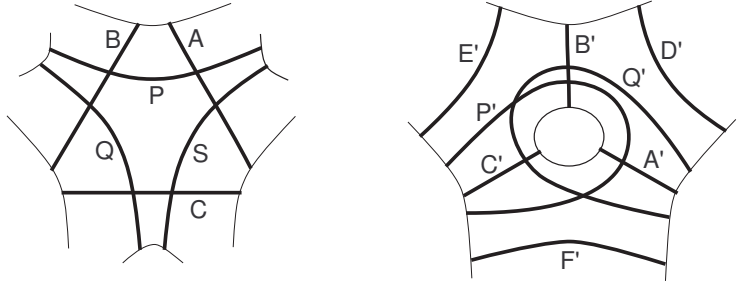}
\caption{Embedded triangle I}
\label{fig100}
\end{center}
\end{figure}

\begin{prop} Let $\{a,b,c\}$ be a $2$-simplex of $\mathcal{A}(R)$. If $\{a,b,c\}$ corresponds to an embedded triangle on $R$, then $\{\lambda(a),\lambda(b),\lambda(c)\}$ corresponds to an embedded triangle on $R$ (see Figure 5).
\label{prop:embedded} \end{prop}

\begin{proof} Let $\Delta$ be an embedded triangle on $R$ with sides corresponding to $A$, $B$, and $C$ representing $a$, $b$ and $c$. Let $T_0 = \{A,B,C\}$. It follows from Proposition \ref{prop:fourtriangles} that we can complete the system of arcs $T_0$ on $R$ to a triangulation $T_1$ of $R$ such that if $\Delta_A$ is the unique triangle of $T_1$ on $R$ different from $\Delta$ having a side corresponding to $A$, $\Delta_B$ is the unique triangle of $T_1$ on $R$ different from $\Delta$ having a side corresponding to $B$, and $\Delta_C$ is the unique triangle of $T_1$ on $R$ different from $\Delta$ having a side corresponding to $C$, then $\Delta$, $\Delta_A$, $\Delta_B$, and $\Delta_C$ are four distinct triangles of $T_1$ on $R$.


Note that $\partial \Delta_A$ is equal to a union of arcs, $A_1$, $X_1$, $B_1$, $Y_1$, $C_1$, and $Z_1$, where $A_1$, $B_1$, and $C_1$ correspond to elements of $T_1$, and each of $X_1$, $Y_1$, and $Z_1$ corresponds to an arc in $\partial R$ or a component of $\partial R$. Without loss of generality, we assume that $A_1$ corresponds to $A$, and $Y_1$ is disjoint from $A_1$.
Similarly, $\partial \Delta_B$ is equal to a union of arcs, $A_2$, $X_2$, $B_2$, $Y_2$, $C_2$, and $Z_2$, where $A_2$, $B_2$, and $C_2$ correspond to elements of $T_1$, and each of $X_2$, $Y_2$, and $Z_2$ corresponds to an arc in $\partial R$ or a component of $\partial R$. Without loss of generality, we assume that $B_2$ corresponds to $B$, and $Z_2$ is disjoint from $B_2$.
Likewise, $\partial \Delta_C$ is equal to a union of arcs, $A_3$, $X_3$, $B_3$, $Y_3$, $C_3$, and $Z_3$, where $A_3$, $B_3$, and $C_3$ correspond to elements of $T_1$, and each of $X_3$, $Y_3$, and $Z_3$ corresponds to an arc in $\partial R$ or a component of $\partial R$. Without loss of generality, we assume that $C_3$ corresponds to $C$, and $X_3$ is disjoint from $C_3$.




Let $P$ be a properly embedded essential arc connecting $Y_1$ to $Z_2$ and intersecting only $A$ and $B$ essentially once and disjoint from other
elements of $T$ as shown in Figure 5. Let $Q$ be a properly embedded essential arc connecting $X_3$ to $Z_2$ and intersecting only $B$ and $C$ essentially once and disjoint from other elements of $T$ as shown in Figure 5. Let $S$ be a properly embedded essential arc connecting $Y_1$ to $X_3$ and intersecting only $A$ and $C$ essentially once and disjoint from other elements of $T$ as shown in Figure 5.

Let $T'$ be a triangulation on $R$ such that $\lambda([T]) = [T']$. Let $A'$, $B'$, and $C'$ be, respectively, the unique representatives of $\lambda(a)$, $\lambda(b)$, and  $\lambda(c)$ in $T'$. Since $\lambda$ is injective $A'$, $B'$, and $C'$ are distinct, and, hence, disjoint and nonisotopic.

Let $p$, $q$, and $s$ be the vertices of $\mathcal{A}(R)$ which are represented by the essential arcs $P$, $Q$, and $S$ on $R$. We see that $\{p, q, s\}$ corresponds to an embedded triangle on $R$.

Let $\sigma$ be the simplex of $\mathcal{A}(R)$ corresponding to the triangulation $T$.
Note that $i(p,x) = 0$ for every vertex $x$ of $\sigma$ other than $a$ and $b$, $i(p,a) = 1$, and $i(p,b) = 1$. Since $\lambda$ is an injective simplicial map, it follows from Proposition \ref{prop:invgeomint} that $i(\lambda(p),y) = 0$ for every vertex $y$ of $\lambda(\sigma)$ other than $\lambda(a)$ and $\lambda(b)$; $i(\lambda(p),\lambda(a)) = 1$, and $i(\lambda(p),\lambda(b)) = 1$. Hence, there exists an arc $P'$ on $R$ representing $\lambda(p)$ such that $P'$ intersects $A'$ and $B'$ essentially once and is disjoint from the other elements of $T'$. Likewise, there exists an arc $Q'$ on $R$ representing $\lambda(q)$ such that $Q'$ intersects $B'$ and $C'$ essentially once and is disjoint from the other elements of $T'$; and there exists an arc $S'$ on $R$ representing $\lambda(s)$ such that $S'$ intersects $C'$ and $A'$ essentially once and is disjoint from the other elements of $T'$.

Since the essential arc $P'$ on $R$ intersects $A'$ and $B'$ essentially once and is disjoint from the other elements of the triangulation $T'$ of $R$, there exists a triangle $\Delta_1$ of $T'$ on $R$ having sides corresponding to $A'$ and $B'$. Similarly, there exists a triangle $\Delta_2$ of $T'$ on $R$ having sides corresponding to $B'$ and $C'$, and a triangle $\Delta_3$ of $T'$ on $R$ having sides corresponding to $C'$ and $A'$. Let the third side of $\Delta_1$ correspond to the element $D'$ of $T'$; the third side of $\Delta_2$ correspond to the element $E'$ of $T'$; and the third side of $\Delta_3$ correspond to the element $F'$ of $T'$.

Suppose, on the one hand, that $D' = C'$. Then $\Delta_1$ is a triangle of $T'$ on $R$ with sides corresponding to $A'$, $B'$ and $C'$. So, $\{\lambda(a),\lambda(b),\lambda(c)\}$ corresponds to an embedded triangle on $R$. Thus, if $D' = C'$, we are done. Likewise, if $E' = A'$ or $F' = B'$, then we are done.

Hence, we may assume that $D' \neq C'$, $E' \neq A'$, and $F' \neq B'$.
Note that, since $A'$, $B'$, and $C'$ are distinct arcs on $R$, $\Delta_1$ has no side corresponding to $C'$. Since $\Delta_2$ has a side corresponding to $C'$, it follows that $\Delta_1 \neq \Delta_2$. Likewise, $\Delta_2 \neq \Delta_3$ and $\Delta_3 \neq \Delta_1$. Hence, $\Delta_1$, $\Delta_2$, and $\Delta_3$ are three distinct components of $R_{T'}$. Since $P$, $Q$, and $S$ are disjoint, $i(p,q) = i(q,s) = i(s,p) = 0$. Since $\lambda$ is a simplicial map, it follows that $i(\lambda(p),\lambda(q)) = i(\lambda(q),\lambda(s)) = i(\lambda(s),\lambda(p)) = 0$. Hence, we may assume that $P'$, $Q'$, and $S'$ are disjoint arcs on $R$. There are three cases to consider, depending on the placement
of the arcs corresponding to $C'$, $E'$, and $F'$ on $\partial \Delta_2$ and $\partial \Delta_3$. These cases are shown in
Figures \ref{fig100} and \ref{fig101}.

Case (i): Assume $A', B', C', D', E', F'$ are as shown in Figure \ref{fig100}. Note that the arc $P'$ on $R$ representing $\lambda(p)$ intersects $B'$ and $A'$ once essentially and is disjoint from $E'$, $D'$, $F'$, and $C'$; and the arc $Q'$ on $R$ representing $\lambda(q)$ intersects $B'$ and $C'$ once essentially and is disjoint from $E'$, $D'$, $F'$, and $A'$. But then we see that $P'$ and $Q'$ intersect essentially (see Figure \ref{fig100}),
which gives a contradiction, since $i(\lambda(p), \lambda(q)) = 0$.

\begin{figure}[htb]
\begin{center}
\epsfxsize=11cm
\epsfbox{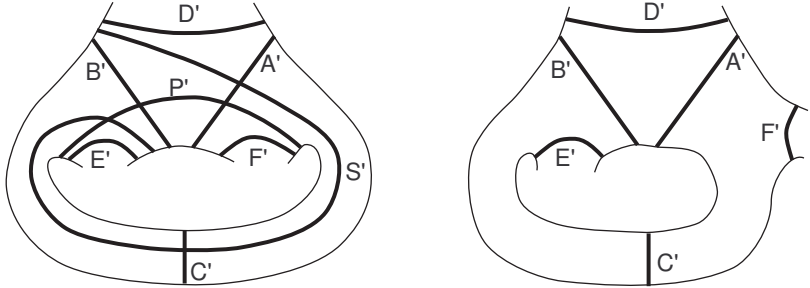}
\caption{Embedded triangle II}
\label{fig101}
\end{center}
\end{figure}

Case (ii): Assume $A', B', C', D', E', F'$ are as shown in the first part of Figure \ref{fig101}. As before, it follows that the arc $P'$ on $R$ representing $\lambda(p)$ intersects $B'$ and $A'$
once essentially and is disjoint from $E'$, $D'$, $F'$, and $C'$; and the arc $S'$ on $R$ representing $\lambda(s)$ intersects $A'$ and $C'$
once essentially and is disjoint from $E'$, $D'$, $F'$, $B'$. But then we see that
$P'$ and $S'$ intersect essentially (see Figure \ref{fig101}),
which gives a contradiction, since $i(\lambda(p),\lambda(s)) = 0$.

The proof for the third case is similar to the proof for Case (ii), (see the second part of
Figure \ref{fig101}).

Hence, we see that either $D' = C'$ or $E' = A'$ or $F' = B'$ and, hence, as argued above, we are done.
\end{proof}

\begin{figure}[htb]
\begin{center}
\epsfxsize=10cm
\epsfbox{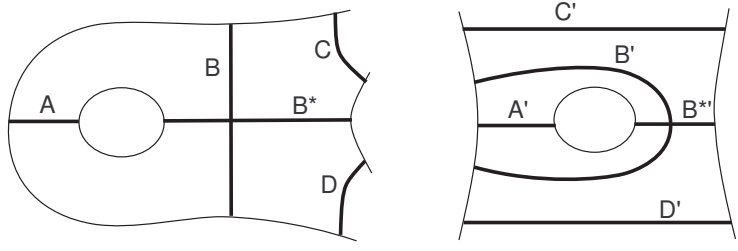}
\caption{Arc configurations I}
\label{fig102}
\end{center}
\end{figure}

Let $(a,b)$ be an ordered $1$-simplex of $\mathcal{A}(R)$. We say that $(a,b)$ {\it corresponds to a non-embedded triangle on $R$} if there exists a non-embedded triangle $\Delta$ on $R$ with sides corresponding to $A$, $B$, and $A$, where $A$ and $B$ represent $a$ and $b$ and $A$ joins two different components of $\partial R$. \label{defn:simpnonembtri}

\begin{prop} Let $(a,b)$ be an oriented edge of $\mathcal{A}(R)$. If $(a,b)$ corresponds to a non-embedded triangle on $R$, then $(\lambda(a),\lambda(b))$ corresponds to a non-embedded triangle on $R$.
\label{prop:nonembedded} \end{prop}

\begin{proof} Let $\Delta$ be a non-embedded triangle on $R$ with sides corresponding to $A$, $B$, where $A$ and $B$ represent $a$ and $b$, and $A$ joins two different components of $\partial R$ (see Figure 7).

Since $R$ is not a pair of pants, there is an embedded triangle $\Delta'$ of $T$ on $R$ having a side corresponding to $B$. Suppose that the other sides of $\Delta'$ correspond to $C$ and $D$ as shown in Figure \ref{fig102}.

Let $B^*$ be as shown in Figure \ref{fig102}, $b^*$ be the vertex of $\mathcal{A}(R)$ corresponding to $B^*$, $c$ be the vertex of $\mathcal{A}(R)$ corresponding to $C$, and $d$ be the vertex of $\mathcal{A}(R)$ corresponding to $D$. Let $\Delta_1$ be the embedded trinangle with sides $A$, $B^*$, and $C$, and $\Delta_2$ be the embedded triangle with sides $A$, $B^*$, and $D$. It follows from Proposition \ref{prop:embedded} that there are embedded triangles $\Delta'_1$ and $\Delta'_2$ on $R$ such that $\Delta'_1$ has sides $A'$, ${B^*}'$, and $C'$, and $\Delta'_2$ has sides $A'$, ${B^*}'$, and $D'$, where $A'$, ${B^*}'$, $C'$, and $D'$ represent $\lambda(a)$, $\lambda(b^*)$, $\lambda(c)$, and $\lambda(d)$.

Since $\lambda$ is an injective simplicial map and $i(b, b^*)=1$, it follows from Proposition \ref{prop:invgeomint} that $i(\lambda(b), \lambda(b^*))=1$. We see that the arc $B'$ representing $\lambda(b)$ can be chosen so that it is disjoint from $A'$,
$C'$, and $D'$ and intersects ${B^*}'$ once. But then, this implies that $A'$ and $B'$ are the sides of a non-embedded triangle on $R$, and $A'$ connects two different components of $\partial R$, (see Figure \ref{fig102}). Since $A'$ and $B'$ represent $\lambda(a)$ and $\lambda(b)$, it follows that $(\lambda(a),\lambda(b))$ corresponds to a non-embedded triangle on $R$.
\end{proof}

\begin{figure}
\begin{center}
\epsfxsize=1.5in \epsfbox{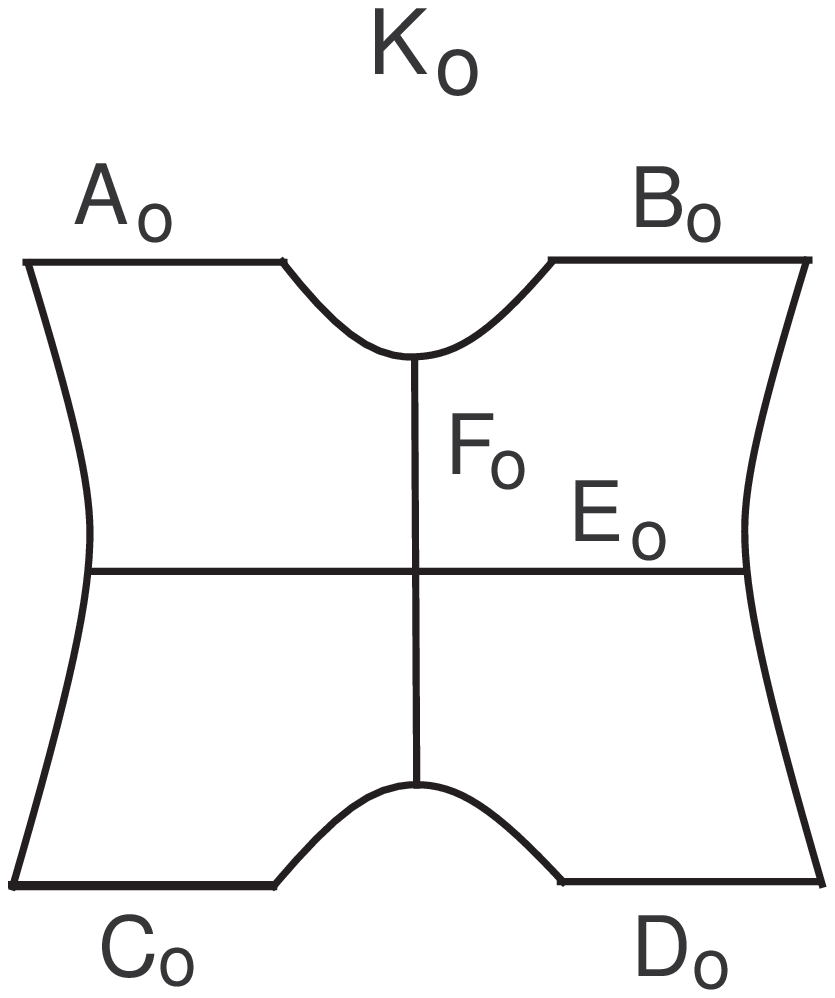} \hspace{0.3cm}
\epsfxsize=1.5in \epsfbox{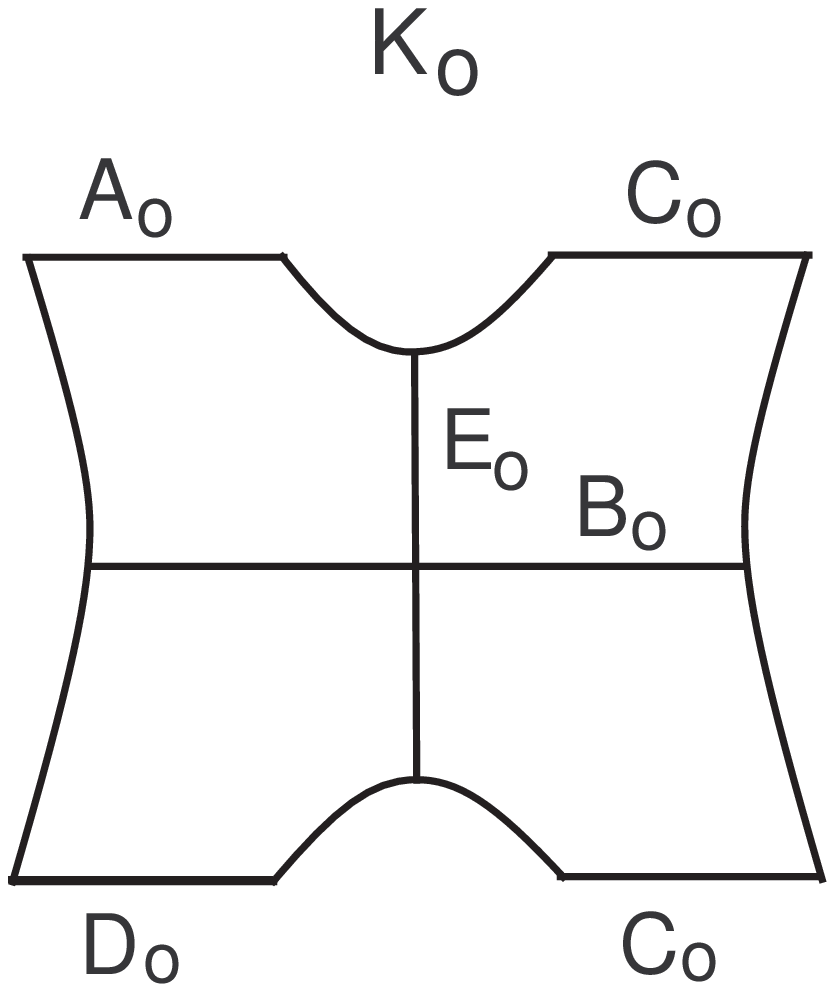} \hspace{0.3cm}
\epsfxsize=1.5in \epsfbox{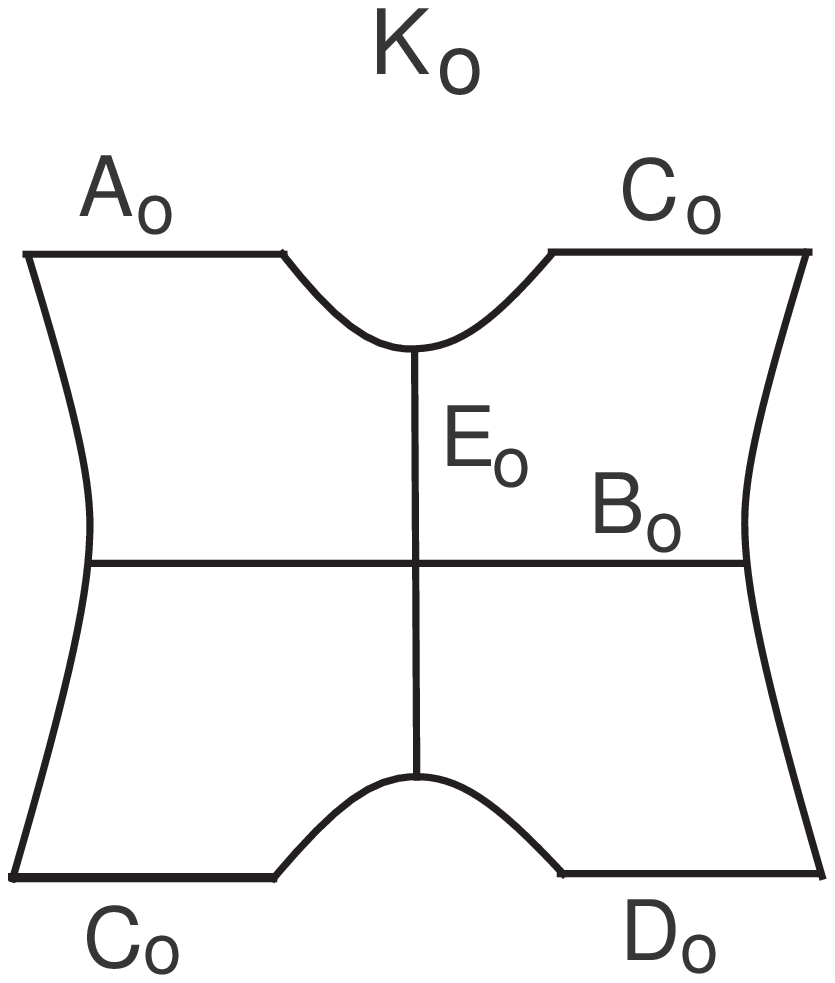}

\hspace{0.cm}  (i)   \hspace{3.7cm}   (ii)   \hspace{3.5cm}    (iii)
\caption{Arc configurations II}
\label{new-fig-1.eps}
\end{center}
\end{figure}

\begin{prop}
\label{9} Let $A, B, C, D$ and $E$ be essential, pairwise disjoint, nonisotopic, properly embedded arcs on $R$. Suppose that there exists a subsurface $K$ of $R$ and a homeomorphism $\phi: (K, A, B, C, D, E) \rightarrow (K_0, A_0, B_0, C_0, D_0, E_0)$ where $K_0$ and $A_0, B_0, C_0, D_0, E_0$ are as shown in Figure \ref{new-fig-1.eps} (i). There exist $A' \in \lambda(a), B' \in \lambda(b), C' \in \lambda(c), D' \in \lambda(d), E' \in \lambda(e)$, $K' \subseteq R$ and a homeomorphism $\chi: (K', A', B', C', D', E') \rightarrow $ $ (K_0, A_0, B_0, C_0, D_0, E_0)$.
\end{prop}

\begin{proof} Let $F_0$ be as shown in Figure \ref{new-fig-1.eps} (i). Let $F = \phi^{-1}(F_0)$. We see that $F$ is an essential properly embedded arc on $R$ such that $F$ intersects $E$ once and $F$ is disjoint from each of $A, B, C, D$. Since $A, B, C, D, E$ are pairwise disjoint, $A, B, E$ and $C, D, E$ form embedded triangles, by using that $\lambda$ is injective and the results of Proposition \ref{prop:embedded}, we can choose
$A' \in \lambda(a), B' \in \lambda(b), C' \in \lambda(c), D' \in \lambda(d), E' \in \lambda(e)$ such that $A', B', C', D', E'$ are pairwise disjoint, and $A', B', E'$ and $C', D', E'$ form embedded triangles. Since $E, F$ intersect once and $F$ is disjoint from each of $A, B, C, D$, we can choose $F' \in \lambda(f)$ such that $E', F'$ intersect once and $F'$ is disjoint from each of $A', B', C', D'$. Since $A, C, F$ form an embedded triangle, $A', C', F'$ form an embedded triangle. Then the result of the proposition follows.
\end{proof}

\begin{prop}
\label{9i} Let $A, B, C, D$ be essential pairwise disjoint nonisotopic properly embedded arcs on $R$. Suppose that there exists a subsurface $K$ of $R$ and a homeomorphism $\phi: (K, A, B, C, D) \rightarrow (K_0, A_0, B_0, C_0, D_0)$ where $K_0$ and $A_0, B_0, C_0, D_0$ are as shown in Figure \ref{new-fig-1.eps} (ii). There exist $A' \in \lambda(a), B' \in \lambda(b), C' \in \lambda(c), D' \in \lambda(d)$, $K' \subseteq R$ and a homeomorphism $\chi: (K', A', B', C', D') \rightarrow (K_0, A_0, B_0, C_0, D_0)$.
\end{prop}

\begin{proof} Let $E_0$ be as shown in Figure \ref{new-fig-1.eps} (ii). Let $E = \phi^{-1}(E_0)$. We see that $E$ is an essential properly embedded arc on $R$ such that $E$ intersects $B$ once and $E$ is disjoint from each of $A, C, D$. Since $A, C, D, E$ are disjoint, $E, C$ form a non-embedded triangle, $A, D, E$ form an embedded triangle, by using that $\lambda$ is injective and the results of Proposition \ref{prop:embedded} and Proposition \ref{prop:nonembedded}, we can choose
$A' \in \lambda(a), E' \in \lambda(e), C' \in \lambda(c), D' \in \lambda(d)$ such that $A', C', D', E'$ are pairwise disjoint, $E', C'$ form a non-embedded triangle and $A', D', E'$ form an embedded triangle. Since $E, B$ intersect once and $B$ is disjoint from each of $A, C, D$, we can choose $B' \in \lambda(b)$ such that $E', B'$ intersect once and $B'$ is disjoint from each of $A', C', D'$, then the result of the proposition follows.
\end{proof}

\begin{prop}
\label{10i} Let $A, B, C, D$ be essential pairwise disjoint nonisotopic properly embedded arcs on $R$. Suppose that there exists a subsurface $K$ of $R$ and a homeomorphism $\phi:
(K, A, B, C, D) \rightarrow (K_0, A_0, B_0, C_0, D_0)$ where $K_0$ and $A_0, B_0, C_0, D_0$ are as shown in Figure \ref{new-fig-1.eps} (iii). There exist $A' \in \lambda(a), B' \in \lambda(b), C' \in \lambda(c), D' \in \lambda(d)$, $K' \subseteq R$ and a homeomorphism $\chi: (K', A', B', C', D') \rightarrow (K_0, A_0, B_0, C_0, D_0)$.
\end{prop}

\begin{proof} Let $E_0$ be as shown in Figure \ref{new-fig-1.eps} (iii). Let $E = \phi^{-1}(E_0)$. Since $A, B, C, D$ are pairwise disjoint, $A, B, C$ form an embedded triangle, and $B, C, D$ form an embedded triangle, by using that $\lambda$ is injective and the results of Proposition \ref{prop:embedded} we choose $A' \in \lambda(a), B' \in \lambda(b), C' \in \lambda(c), D' \in \lambda(d)$ such that $A', B', C', D'$ are pairwise disjoint, $A', B', C'$ form an embedded triangle and $B', C', D'$ form an embedded triangle. Since $E, B$ intersect once and $E$ is disjoint from each of $A, C, D$, we can choose $E' \in \lambda(e)$ such that $E', B'$ intersect once and $E'$ is disjoint from each of $A', C', D'$. Since $A, C, E$ form an embedded triangle, $A', C', E'$ form an embedded triangle. Then the result of the proposition follows.
\end{proof}

\begin{prop} Let $\sigma$ be a simplex of $\mathcal{A}(R)$ corresponding to a triangulation $T$ of $R$. Let $T'$ be a triangulation of $R$ such that $\lambda(\sigma)$ is the simplex of $\mathcal{A}(R)$ corresponding to $T'$. For each arc $J$ of $T$ let $J'$ be the unique arc of $T'$ such that $\lambda([J])= [J']$. Then there exist a homeomorphism $H : R \rightarrow R$ such that $H(J) = J'$ for each arc $J$ of $T$. \label{prop:ordtopequiv} \end{prop}

\begin{proof} Let $R_1$ denote the surface obtained from $R$ by cutting $R$ along $T$ and $R_2$ be the surface obtained from $R$ by cutting $R$ along $T'$. $R_1$ and $R_2$ both have $N$ components, where $N = 4g + 2b - 4$. Let $\{\Delta_i | 1 \leq i \leq N\}$ be the $N$ distinct components of $R_1$.

Since $R$ is not a pair of pants or a torus with two holes, no two distinct components of $R_1$ can have sides corresponding to the same elements of $T$. Likewise, no two distinct components of $R_2$ can have sides corresponding to the same elements of $T'$.

Let $i$ be an integer with $1 \leq i \leq N$. Since $\Delta_i$ is a component of $R_1$, $\Delta_i$ is a triangle of $T$ with sides corresponding to elements $A_i$, $B_i$, and $C_i$ of $T$. Let $a_i$, $b_i$, and $c_i$ be the vertices of $\mathcal{A}(R)$ represented by $A_i$, $B_i$ and $C_i$. Then, $\{a_i,b_i,c_i\}$ corresponds to a triangle on $R$. Hence, by Propositions \ref{prop:embedded} and \ref{prop:nonembedded} , $\{a'_i,b'_i,c'_i\}$ corresponds to a triangle on $R$, where $a'_i = \lambda(a_i)$, $b'_i = \lambda(b_i)$, and $c'_i = \lambda(c_i)$. Let $A'_i$, $B'_i$, and $C'_i$ be the unique elements of $T'$ which represent $a'_i$, $b'_i$, and $c'_i$. It follows that there exists a unique triangle $\Delta'_i$ of $T'$ on $R$ with sides corresponding to $A'_i$, $B'_i$, and $C'_i$.

Moreover, the correspondence $\Delta_i \mapsto \Delta'_i$ establishes a bijection from the set of exactly $N$ distinct components $\{\Delta_i | 1 \leq i \leq N\}$ of $R_1$ to the set of exactly $N$ distinct components $\{\Delta'_i | 1 \leq i \leq N\}$ of $R_2$.

Suppose, on the one hand, that $\Delta_i$ is embedded. Then, by Proposition \ref{prop:embedded}, $\Delta'_i$ is embedded. Let $J_i$, $K_i$, and $L_i$ be the arcs in $\Delta_i$ corresponding to $A_i$, $B_i$, and $C_i$, and $J'_i$, $K'_i$, and $L'_i$ be the arcs in $\Delta'_i$ corresponding to $A'_i$, $B'_i$, and $C'_i$. Note that there exists a homeomorphism $F_i : (\Delta_i,J_i,K_i,L_i) \rightarrow (\Delta'_i,J'_i,K'_i,L'_i)$ which is well-defined up to relative isotopies. In particular, the orientation type of $F_i : (\Delta_i,J_i,K_i,L_i) \rightarrow (\Delta'_i,J'_i,K'_i,L'_i)$ (i.e. whether it is orientation-reversing or orientation-preserving) is fixed.

Suppose, on the other hand, that $\Delta_i$ is non-embedded. Let $J_i$, $K_i$, and $L_i$ be the arcs in $\Delta_i$ corresponding to $A_i$, $B_i$, and $C_i$, and $J'_i$, $K'_i$, and $L'_i$ be the arcs in $\Delta'_i$ corresponding to $A'_i$, $B'_i$, and $C'_i$. We may assume that $A_i = C_i$, so that $A_i$ joins two different components of $\partial R$, and $B_i$ joins a component of $\partial R$ to itself. Then, by Proposition \ref{prop:nonembedded}, $\Delta'_i$ is non-embedded, $A'_i = C'_i$, $A'_i$ joins two different boundary components of $\partial R$, and $B'_i$ joins a component of $\partial R$ to itself.

In this situation there is an ambiguity in the choice of $J_i$ and $L_i$. After all, $J_i$ and $L_i$ both correspond to $A_i$ (i.e. $J_i$ and $L_i$ both correspond to $C_i$). Likewise, there is an ambiguity in the choice of $J'_i$ and $L'_i$. Suppose that $(J_i,L_i,J'_i,L'_i)$ has been specified. Then there exist a homeomorphism $F_i : (\Delta_i,J_i,K_i,L_i) \rightarrow (\Delta'_i,J'_i,K'_i,L'_i)$ which is well-defined up to relative isotopies, and a homeomorphism $F^*_i :  (\Delta_i,J_i,K_i,L_i) \rightarrow $ $ (\Delta'_i,L'_i,K'_i,J'_i)$ which is well-defined up to isotopies. In particular, in this situation, the orientation type of $F_i : (\Delta_i,J_i,K_i,L_i) \rightarrow (\Delta'_i,J'_i,K'_i,L'_i)$ is fixed; the orientation type of $F^*_i : (\Delta_i,J_i,K_i,L_i) \rightarrow (\Delta'_i,L'_i,K'_i,J'_i)$ is fixed; and $F_i : (\Delta_i,J_i,K_i,L_i) \rightarrow (\Delta'_i,J'_i,K'_i,L'_i)$  and $F^*_i : (\Delta_i,J_i,K_i,L_i) \rightarrow (\Delta'_i,L'_i,K'_i,J'_i)$ have opposite orientation types.

By using Propositions \ref{prop:embedded}, \ref{prop:nonembedded}, \ref{9}, \ref{9i} and \ref{10i}, we can choose homeomorphisms $G_i : \Delta_i \rightarrow \Delta'_i$, $1 \leq i \leq N$, where $G_i$ is isotopic to $F_i$, if $\Delta_i$ is embedded, and $G_i$ is isotopic to either $F_i$ or $F^*_i$, if $\Delta_i$ is non-embedded, so that the unique homeomorphism $G : R_1 \rightarrow R_2$ whose restriction to $\Delta_i$ is equal to $G_i$, $1 \leq i \leq N$, covers a homeomorphism $H : R \rightarrow R$. Roughly speaking, these propositions ensure that the homeomorphisms $F_i : \Delta_i \rightarrow \Delta'_i$ and $F_j : \Delta_j \rightarrow \Delta'_j$ associated to embedded triangles $\Delta_i$ and $\Delta_j$ which have sides corresponding to the same element of $T$, can be isotoped by a relative isotopy to agree, relative to the natural quotient maps, $q_1: R_1 \rightarrow R$ and $q_2 : R_2 \rightarrow R$. In other words, the restrictions of $F_i$ and $F_j$ to pairs of sides which correspond to the same element of $T$, which restrictions may be identified, via $q_1$ and $q_2$, to homeomorphisms from a fixed element of $T$ to a fixed element of $T'$, have the same orientation type as such homeomorphisms between fixed elements of $T$ and $T'$.

When $\Delta_i$ is nonembedded, this condition on compatibility of orientation types of restrictions on pairs of sides which correspond to the same element of $T$ can be realized on all such pairs by making the appropriate choice of either $F_i$ or $F^*_i$, $1 \leq i \leq N$.

Once the correct choices are made so that this compatibility of orientations is realized, we may isotope the chosen homeomorphisms, $F_i$ or $F^*_i$, to homeomorphisms $G_i$ which agree, as homeomorphisms between fixed elements of $T$ and $T'$, on all pairs of sides which correspond to the same element of $T$.

It follows that $H : R \rightarrow R$ is a homeomorphism which maps each element $J$ of $T$ to the corresponding element $J'$ of $T'$, completing the proof.\end{proof}\\

We will need the following strong form of connectivity for $\mathcal{A}(R)$ stated as the ``Connectivity Theorem for Elementary Moves'' in Mosher \cite{Mos}. See also Corollary 5.5.B in Ivanov's survey article on Mapping Class Groups \cite{ivanov2}.

\begin{theorem}(Connectivity Theorem for Elementary Moves, \cite{Mos}) Suppose that $R$ is not a disc or an annulus. Then any two triangulations of $R$ are related by a finite sequence of elementary moves. More precisely, if $T$ and $T'$ are triangulations of $R$ and $\sigma$ and $\sigma'$ are the simplices of $\mathcal{A}(R)$ corresponding to $T$ and $T'$, then there exists a sequence of simplices $\sigma_i$, $1 \leq i \leq N$, such that $\sigma_1 = \sigma$, $\sigma_N = \sigma'$, and for each integer $i$ with $1 \leq i < N$, $\sigma_{i+1}$ is obtained from $\sigma_i$ by an elementary move. \label{thm:connectivity} \end{theorem}

\begin{rem} The statement of Theorem \ref{thm:connectivity} in Mosher \cite{Mos} is in terms of ideal triangulations of a punctured surface $(S,P)$ rather than triangulations of $R$. For our purposes here, we let $S$ be the closed surface of genus $g$ obtained from $R$ by attaching a disc $D_i$ to each component $\partial_i$ of $\partial R$, $1 \leq i \leq b$, and $P$ be a set of points, $x_i$, $1 \leq i \leq p$, with $x_i$ in the interior of $D_i$, $1 \leq i \leq b$. Then we may relate triangulations of $R$ as defined in this paper to ideal triangulations of $(S,P)$ as defined in Mosher \cite{Mos} by ``coning off'' arcs on $R$ to arcs or loops on $S$ joining points in $P$ to points in $P$. In this way, we obtain the above restatement of the Connectivity Theorem for Elementary Moves in a form suitable for our purposes in this paper.\label{rem:connectivity2} \end{rem}

\begin{theorem} Suppose that $(g, b) \neq (0,1), (g, b) \neq (0,2)$, $(g, b) \neq (0,3)$ and $(g, b) \neq (1,1)$. Let $\lambda: \mathcal{A}(R) \rightarrow \mathcal{A}(R)$ be an injective simplicial map. Then $\lambda : \mathcal{A}(R) \rightarrow \mathcal{A}(R)$ is geometric (i.e. there exists a homeomorphism $H : R \rightarrow R$ such that for every essential arc $A$ on $R$, $\lambda([A]) = [H(A)]$).\label{thm:injectivegeometric} \end{theorem}

\begin{proof} Let $\sigma$ be a maximal simplex of $\mathcal{A}(R)$. Let $T$, $T'$, and $H : R \rightarrow R$ be as in Proposition \ref{prop:ordtopequiv}. Let $\psi = H_*^{-1} \circ \lambda : \mathcal{A}(R) \rightarrow \mathcal{A}(R)$. We see that
$\psi(x) = x$ for each vertex $x$ of $\sigma$. Recall that (i) each vertex of $\mathcal{A}(R)$ is contained in a codimension zero face of $\mathcal{A}(R)$, (ii) each codimension one face of $\mathcal{A}(R)$ is contained in one or two codimension zero faces of $\mathcal{A}(R)$, and (iii) Theorem \ref{thm:connectivity}, Mosher's ``Connectivity by Elementary Moves'' holds. It follows from these facts that $\psi = id_{\mathcal{A}(R)} : \mathcal{A}(R) \rightarrow \mathcal{A}(R)$. Hence, $\lambda = H_* : \mathcal{A}(R) \rightarrow \mathcal{A}(R)$. That is to say, $\lambda$ is geometric, being induced by the self-homeomorphism $H : R \rightarrow R$.
\end{proof}


\begin{theorem} Suppose that $(g, b) \neq (0,1), (g, b) \neq (0,2), (g, b) \neq (0,3),$ and $(g, b) \neq (1,1)$. Then $Aut(\mathcal{A}(R))$ is naturally isomorphic to the extended mapping class group $\Gamma^*(R)$.\end{theorem}

\begin{proof} By Theorem \ref{thm:injectivegeometric}, the natural representation $\rho : \Gamma^*(R) \rightarrow Aut(\mathcal{A}(R))$ is surjective. Let $h$ be an element of $ker(\rho)$ and $H : R \rightarrow R$ be a homeomorphism of $R$ representing $h$.
$H$ preserves the isotopy class of every essential arc on $R$. Since $R$ is not a disc or an annulus, there exists a triangulation $T$ of $R$. Since $H: R \rightarrow R$ preserves the isotopy class of every essential arc on $R$, we may isotope $H : R \rightarrow R$ to a homeomorphism $H_0 : R \rightarrow R$ such that, for each element $J$ of $T$, $H_0(J) = J$.

Since $R$ is not a pair of pants, there exists an embedded triangle $\Delta$ of $T$ on $R$ with sides $A$, $B$, and $C$. Since $A$, $B$ and $C$ are elements of $T$, there exists an embedded triangle $\Delta'$ of $T$ on $R$ with sides corresponding to $H_0(A)$, $H_0(B)$, and $H_0(C)$.

Since $A$, $B$, and $C$ are elements of $T$, $H_0(A) =  A$, $H_0(B) =  B$, and $H_0(C) = C$. It follows that $\Delta'$ and $\Delta$ are triangles on $R$ with sides $A$, $B$ and $C$. Since $R$ is not a pair of pants $\Delta' = \Delta$. So, $H_0(\Delta)$ is equal to $\Delta$ and, hence, the homeomorphism $H_0 : R \rightarrow R$ restricts to a homeomorphism $H_0| : (\Delta, A, B, C) \rightarrow (\Delta, A, B, C)$. We may isotope $H_0 : R \rightarrow R$ relative to the union $|T|$ of the elements of $T$, to a homeomorphism $H_1 : R \rightarrow R$ which restricts to the identity map $H_1| = id_{\Delta} : \Delta \rightarrow \Delta$ of $\Delta$.

Note that any other triangle $\Delta''$ of the triangulation $T$ of $R$ is connected to the triangle $\Delta$ of the triangulation $T$ of $R$ by a sequence of triangles which have sides corresponding to the same element of $T$. Since $H : R \rightarrow R$ is orientation-preserving, it follows, by a finite induction argument, that we may construct a sequence of homeomorphisms, $H_i : R \rightarrow R$, $0 \leq i \leq N$, with $N$ equal to the number of triangles of $T$ on $R$, such that $H_0$ preserves each element of $T$ and is isotopic on $R$ to $H$; $H_1$ preserves each element of $T$, fixes each point of at least one triangle of $T$ on $R$, and is isotopic on $R$ to $H_0$ relative to $|T|$; and for each integer $i$ with $2 \leq i \leq N$, $H_i$ preserves each element of $T$, fixes each point of at least $i$ triangles of $T$, and is isotopic on $R$ to $H_{i - 1}$ relative to the union of $|T|$ with $i-1$ triangles of $T$ fixed pointwise by $H_{i-1}$.

Since $N$ is equal to the number of triangles of $T$ on $R$, it follows that $H_N = id_R : R \rightarrow R$. Since $H : R \rightarrow R$ is, by induction, isotopic to $H_N : R \rightarrow R$, it follows that $H : R \rightarrow R$ is isotopic to $id_R : R \rightarrow R$. Hence, we have $\Gamma^*(R) \cong Aut(\mathcal{A}(R))$.\end{proof}

\noindent Elmas Irmak, Department of Mathematics and Statistics, Bowling Green State University, Bowling Green, OH 43403, USA,
eirmak@bgsu.edu.\\

\noindent John D. McCarthy, Department of Mathematics, Michigan State University, East Lansing, MI 48824-1027, USA, mccarthy@math.msu.edu.

\end{document}